\documentclass[twocolumn]{autart}    

\usepackage[dvips]{epsfig}    
\usepackage{mathrsfs}

\usepackage{bbm}
\usepackage{color}
\usepackage[comma,numbers,square,sort&compress]{natbib}
\usepackage{amsmath,amssymb,amsfonts}
\usepackage{epstopdf}
\usepackage{mathtools}

\usepackage{tikz}
\usepackage{verbatim}
\usepackage{graphicx}
\usepackage{subfigure}
\usepackage{graphics} 
\usepackage{epsfig} 

\usepackage{marvosym}

\usetikzlibrary{shapes,arrows}

\newtheorem{theorem}{Theorem}

\newtheorem{lemma}{Lemma}
\newtheorem{defin}{Definition}
\newtheorem{examp}{Example}

\newcommand\Set[2]{%
 \left\{#1\ \middle\vert\ #2 \right\}}

\begin{document}

\begin{frontmatter}

\title{Controllability of Heterogeneous Directed Networked MIMO Systems\thanksref{footnoteinfo}}

\thanks[footnoteinfo]{A preliminary version of this work was presented at IWCSN 2017 \cite{wang2017con}. This work was supported in part by the National Natural Science Foundation of China under Grant 61473240 and in part by the Natural Science Foundation of Fujian Province under Grant 2017J01119, as well as by the Hong Kong Research Grants Council under the GRF Grant CityU 11200317.}
\thanks[correspondingauthor]{Correspondence should be addressed to Peiru~Wang (wangpeiru.hfut@gmail.com).}

\author[school]{Linying~Xiang}\ead{xianglinying@neuq.edu.cn},    
\author[Xiamen]{Peiru~Wang}\ead{wangpeiru.hfut@gmail.com},               
\author[school,keylab]{Fei~Chen}\ead{fei.chen@ieee.org},  
\author[Hongkong]{Guanrong~Chen}\ead{eegchen@cityu.edu.hk}

\address[school]{School of Control Engineering, Northeastern University at Qinhuangdao,
Qinhuangdao, 066004, China}  
\address[Xiamen]{Department of Automation, Xiamen University, Xiamen, 361005, China}             
\address[keylab]{State Key Laboratory of Synthetical Automation for
Process Industries, Northeastern University, Shenyang, 110004, China}        
\address[Hongkong]{Department of Electronic Engineering, City University of Hong Kong, 83 Tat Chee Avenue, Kowloon, Hong Kong SAR, China}

\begin{keyword}                           
Controllability; heterogeneous network; MIMO system; directed graph.              
\end{keyword}                             

\begin{abstract}                          
This paper studies the controllability of networked multi-input-multi-output (MIMO) systems, in which the network topology is weighted and directed, and the nodes are heterogeneous higher-dimensional linear time-invariant (LTI) dynamical systems. The primary objective is to search for controllability criteria beyond those already known for homogeneous networks. The focus is on the effects of the network topology, node dynamics, external control inputs, as well as the inner interactions on the network controllability. It is found that a network of heterogeneous systems can be controllable even if the corresponding homogeneous network topology is uncontrollable. The finding thus unravels another fundamental property that affects the network controllability---the heterogeneity of the node dynamics.
A necessary and sufficient condition is derived for the controllability of heterogeneous networked MIMO LTI systems. For some typical cases, necessary and/or sufficient controllability conditions are specified and presented on the node dynamics, inner interactions, as well as the network topology.
\end{abstract}

\end{frontmatter}

\section{Introduction}
The ultimate goal of understanding complex networked systems is to control their functioning and behavior. To fully control a network of dynamical systems, one should first determine whether or not the network is controllable \cite{liu2011}, \cite{tanner2004}. Controllability, one of the fundamental concepts in control theory, quantifies the ability to steer a dynamical system from any initial state to any final state within finite time \cite{kalman1963}. The classical notion of controllability in control theory mainly emphasizes on the inherent dynamics of a single higher-dimensional system, which is pertinent to its microcosmic component dynamics. In the big-data era and omni-networking world today, the traditional control theory is encountered more and more with large-scale networked systems \cite{chen2017}, where nodes are the elements of the network and edges represent the interactions among them. Typical examples include the WWW, the Internet, transportation networks, wireless communication networks, power grids, social networks, and biological networks, to name just a few. Perturbations to one node in the network can regulate the states of the other nodes through their local interactions. This property enables the possibility of controlling the whole network by manipulating the states of only a subset of the nodes, for which the topology of the network is crucial. Therefore, it is of both theoretical and practical importance to explore the controllability of a complex network of dynamical systems, which helps understand, predict and optimize the collective behavior of complex networked systems from the macrocosmic perspective.

Consider a homogeneous network of $N$ LTI systems (nodes) described by $\dot{x}=Ax+Bu$, where $x\in \mathbb{R}^{N}$ is the state vector of the nodes, $u\in \mathbb{R}^{M}$ is the input vector, $A\in \mathbb{R}^{N\times N}$ is the adjacency matrix of the underlying network, and $B\in \mathbb{R}^{N\times M}$ is the input matrix identifying the nodes that are directly under control.  It is convenient to examine the controllability via Kalman's rank criterion \cite{kalman1963} when the network size is small: a node system is (state) controllable if and only if the controllability matrix $Q=[B,AB,\dots,A^{N-1}B]$ has full
row-rank. Yet, this routine generally fails to work for large-scale networked systems due to high computational cost. Additionally, to numerically check the rank condition, one has to know the precise values of the parameters in matrices $A$ and $B$, which are seldom possible in practice. For this reason, structural controllability was proposed in \cite{lin1974} to relax the above limitation. The concept highlights the role of the underlying network structure in controllability, where the system parameters can be either fixed zeros or independent nonzero parameters. Under the framework of the structural controllability theory, one can determine the network controllability even if the exact values of the edge weights are not available. Based on the structural controllability theory, the controllability of large-scale weighted and directed single-input-single-output (SISO) networks was investigated in \cite{liu2011} and a  minimum inputs theorem was established to identify the minimum number of driver nodes that need to be controlled by external signals to ensure the network controllability. The basic idea is to seek the unmatched nodes in the network using a maximum matching algorithm \cite{hopcroft1973}. Since then, the issue of network controllability for complex dynamical systems has become a focal subject in network science \cite{yao2017}, and numerous works have been reported from rather diverse perspectives on such topics as control capacity \cite{jia2013}, edge dynamics \cite{nepusz2012}, \cite{pang2017}, optimization \cite{wang2012}, \cite{xiao2014}, \cite{ding2016}, control energy \cite{yan2015energy} \cite{sun2013energy}, \cite{klickstein2017energy}, \cite{sun2017energy},  exact controllability \cite{yuan2013}, \cite{yuan2014}, \cite{li2014} and robustness \cite{nie2014robustness}, \cite{yan2016robustness}, \cite{lou2018robustness}.

Recently, it has also been revealed that node dynamics is another significant factor affecting system controllability in addition to the network topology. In \cite{cowan2012}, it was pointed out that the main results in \cite{liu2011} depend heavily on a critical assumption: each node has an infinite time constant (i.e., each node is treated as a pure integrator); however, the real networks considered therein include food webs, power grids, electronic circuits, regulatory networks, and neuronal networks, which typically have finite time constants. Indeed, by analyzing the structural controllability of directed networks with LTI nodal dynamics, it was found that only a single input is required to ensure the network controllability. In \cite{zhao2015}, the synergistic effect of the network topology was investigated along with the so-called $d$-order individual dynamics on exact controllability. A global symmetry relationship was found, which accounts for the invariance of controllability with respect to exchanging the densities of any two different types of dynamic units, irrespective of the network topology. More general results on MIMO node systems can be found in \cite{wang2016}, \cite{wang2017}, \cite{hao2017} and the recent survey\cite{xiang2019}.

It should be pointed out that the above-reviewed works all assume that all nodes in the networked systems have identical self-dynamics \cite{zhao2015}, \cite{wang2016}, \cite{wang2017}, \cite{hao2017}. However, homogeneity is only an ideal assumption that can rarely be satisfied for practical systems. In a more realistic setting, the heterogeneity of the networked systems cannot be neglected, for which the results of homogeneous networked systems do not apply, therefore new approaches need to be developed \cite{xiang2013}\cite{zhao2017controllability}. Motivated by the above discussions, the controllability of heterogeneous networked MIMO systems is investigated in this paper. It will become clear in Section~\ref{sec:two_examp}
that the controllability of a heterogeneous network differs dramatically from that of a homogeneous
network.

The rest of the paper is organized as follows. Section~II introduces the notation and graph theory to be used throughout the paper. Section~III describes the problem to be investigated and presents two examples, which demonstrate that the controllability of heterogeneous networks is essentially different from that of homogeneous networks. In Section~IV, the main results are presented for the controllability of heterogeneous MIMO systems. In Section~V, heterogeneous networks with controllable node systems are discussed in more detail. Finally, Section~VI concludes the paper.



\section{Notation and graph theory}
\subsection{Notation}
Let $\mathbb{R}$ ($\mathbb{C}$) denote the set of real (complex) numbers, $\mathbb{R}^{n}$ ($\mathbb{C}^{n}$) denote the vector space of the $n$-dimensional real (complex) vectors, and $\mathbb{R}^{n\times m}$ ($\mathbb{C}^{n\times m}$) denote the set of $n\times m$ real (complex)  matrices, with $I_{N}$ being the $N\times N$ identity matrix and $diag\{a_{1},...,a_{N}\}$ being the $N\times N$ diagonal matrix. Let $\sigma(A)$  denote the set of all the eigenvalues of matrix $A$.


\subsection{Graph theory}
A directed graph $G=(V,E)$ consists of a node set $V=\{\nu_{1},\dots,\nu_{N}\}$ and an edge set $E=\{(\nu_{j},\nu_{i})\}$. Let $W=[\omega_{ij}]\in R^{N\times N}$ denote the weighted adjacency matrix of the graph, where $\omega_{ij}\neq 0$ if $(\nu_{j},\nu_{i})\in E$ and $\omega_{ij}=0$ otherwise. In this paper, simple directed and weighted graphs are considered; that is, self-loops and multiple edges are excluded. A chain network consisting of $N$ nodes is a directed path from node $1$ to node $N$.

\section{Problem statement}
\subsection{A heterogenous network model}
Consider the following heterogeneous network model:
\begin{equation}
 \left\{ \begin{array}{ll}
 \dot{x}_i=A_ix_i+\sum_{j=1}^{N}\omega_{ij}Hy_{j}+\delta_{i}B_{i}u_{i},\\
 y_{i}=C_{i}x_{i},~~~~i=1,\dots,N,\\
 \end{array} \right. \label{eq:sys_mod}
\end{equation}
where $x_{i}\in \mathbb{R}^{n}$, $u_{i}\in \mathbb{R}^{p}$  and  $y_{i}\in \mathbb{R}^{m}$ denote, respectively, the state, input and output of node $i$. The matrices $A_i$, $B_i$ and $C_i$ are, respectively, the state, input and output matrices of node $i$.  The matrix $H\in \mathbb{R}^{n\times m}$ describes the coupling among different components. The weighted adjacency matrix $W=[\omega_{ij}]\in \mathbb{R}^{N\times N}$ represents the network topology and defines the strengths of the interactions among the nodes. The binary variable $\delta_i$ indicates whether node $i$ is under control, i.e., $\delta_{i}=1$ if node $i$ is under control, and $\delta_{i}=0$ otherwise.

Define $x=[x_{1}^{T},\dots,x_{N}^{T}]^{T}$ and $u=[u_{1}^{T},\dots,u_{N}^{T}]^{T}$. Additionally, define the matrices $\Phi=\Lambda+\Gamma$, where  $\Lambda=\text{diag}(\{A_{1},\dots,A_{N}\})$ and $\Gamma=[\Gamma_{ij}]$ with $\Gamma_{ij}=\omega_{ij}HC_{j}$,  and $\Psi=\text{diag}(\{\delta_{1}B_{1},\dots,\delta_{N}B_{N}\})$. One can rewrite the network model \eqref{eq:sys_mod} in a corresponding compact system form, as
\begin{equation}
 \dot x=\Phi x+\Psi u. \label{eq:sys_mod_stacked}
\end{equation}

Let $\Delta=\text{diag}(\{\delta_{1},...,\delta_{N}\})$. The following definition characterizes whether a network topology is controllable.
\begin{defin}
   The network topology is said to be controllable if and only if $(W,\Delta)$ is a controllable matrix pair.
\end{defin}

It is noted that the controllability of network topology \eqref{eq:sys_mod} differs from the controllability of \eqref{eq:sys_mod_stacked}. For a homogeneous network, the controllability of the network topology is necessary for the controllability of the corresponding system. However, it might not be necessary for heterogeneous networks as will be seen in the following subsection.

\subsection{Two comparative examples }
\label{sec:two_examp}

In this section, two examples are presented to highlight the difference between the controllabilities of heterogeneous networks and homogeneous networks.

\begin{examp}
Consider a directed chain network of three nodes, shown in Fig.~\ref{fig:chain_netw}. Let $\omega_{21}=\omega_{32}=1$, $\delta_{1}=1$, and $\delta_{2}=\delta_{3}=0$.
It is straightforward to verify that $(W,\Delta)$ is controllable. Suppose that each node has the following identical matrices:
\begin{figure}[htp]
  \centering

\begin{tikzpicture}[-,>=stealth',shorten >=1pt,auto,node distance=1.6cm,thick,
    node/.style={scale=0.8, circle,fill=white!20,draw,font=\sffamily\tiny\bfseries},
    input/.style={scale=0.8, circle,fill=white!20,draw=white,font=\sffamily\tiny\bfseries},
    function/.style={scale=0.8,draw,fill=black!20,font=\sffamily\tiny\bfseries}
    ]
  \node[input] (v1) {$u_{1}$};
  \node[node] (v2) [below of=v1] {$1$} edge [<-] (v1);
  \node[node] (v3) [below of=v2] {$2$}  edge [<-] (v2);
  \node[node] (v4) [below of=v3] {$3$} edge [<-] (v3);

\end{tikzpicture}
  \caption{\textbf{A directed chain network with three nodes. }\label{fig:chain_netw} }

\end{figure}
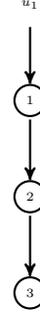

$$A_{1}=A_{2}=A_{3}=
\begin{bmatrix}
  1 & 0 \\
  2 & 1 \\
\end{bmatrix},$$
$$B_{1}=B_{2}=B_{3}=
\begin{bmatrix}
  1 & 2 \\
  0 & 1 \\
\end{bmatrix},$$
$$C_{1}=C_{2}=C_{3}=
\begin{bmatrix}
  1 & 0 \\
  0 & 2 \\
\end{bmatrix},$$
$$H=
\begin{bmatrix}
  1 & 0 \\
  0 & 1 \\
\end{bmatrix}.$$
By checking the rank of $[\Psi,\Phi\Psi,\Phi^{2}\Psi, \Phi^3 \Psi, \Phi^4 \Psi, \Phi^5 \Psi]$, defined in system (2), it is clear that the networked system is controllable.
However, when the node dynamics are heterogeneous, e.g.,
$$A_{1}=
\begin{bmatrix}
  1 & 0 \\
  2 & 1 \\
\end{bmatrix}
,~
A_{2}=
\begin{bmatrix}
  1 & 0 \\
  0 & 1 \\
\end{bmatrix}
,~
A_{3}=
\begin{bmatrix}
  1 & 2 \\
  0 & 1 \\
\end{bmatrix},$$
$$B_{1}=
\begin{bmatrix}
  1 & 0 \\
  0 & 1 \\
\end{bmatrix}
,~
B_{2}=
\begin{bmatrix}
  1 & 0 \\
  0 & 2 \\
\end{bmatrix}
,~
B_{3}=
\begin{bmatrix}
  2 & 0 \\
  0 & 1 \\
\end{bmatrix},$$
$$C_{1}=
\begin{bmatrix}
  0 & 1 \\
  0 & 0 \\
\end{bmatrix}
,~
C_{2}=
\begin{bmatrix}
  0 & 2 \\
  0 & 0 \\
\end{bmatrix}
,~
C_{3}=
\begin{bmatrix}
  0 & 3 \\
  0 & 0 \\
\end{bmatrix},$$
$$H=
\begin{bmatrix}
  1 & 0 \\
  0 & 1 \\
\end{bmatrix},$$
the corresponding system is uncontrollable.
\end{examp}

\begin{examp}
Consider a directed tree of three nodes, shown in Fig.~\ref{fig:tree_netw}. Let $\omega_{21}=\omega_{31}=1$,  $\delta_{1}=1,$ and $\delta_{2}=\delta_{3}=0$.
It is straightforward to verify that $(W,\Delta)$ is uncontrollable. {When the three nodes have the same dynamics, i.e., $A_{1}=A_{2}=A_{3}$, $B_{1}=B_{2}=B_{3}$, $C_{1}=C_{2}=C_{3}$}, the networked system is uncontrollable. However, for the following heterogeneous node dynamics:
\begin{figure}[htp]
\centering
\begin{tikzpicture}[-,>=stealth',shorten >=1pt,auto,node distance=1.6cm,thick,
    node/.style={scale=0.8, circle,fill=white!20,draw,font=\sffamily\tiny\bfseries},
    input/.style={scale=0.8, circle,fill=white!20,draw=white,font=\sffamily\tiny\bfseries},
    function/.style={scale=1.1,draw,fill=black!20,font=\sffamily\tiny\bfseries}
    ]
  \node[input] (v1) {$u_{1}$};
  \node[node] (v2) [below of=v1] {$1$} edge [<-] (v1);
  \node[node] (v3) [below left of=v2] {$2$}  edge [<-] (v2);
  \node[node] (v4) [below right of=v2] {$3$} edge [<-] (v2);
\end{tikzpicture}
\caption{A directed tree with three nodes.  \label{fig:tree_netw}}
\end{figure}
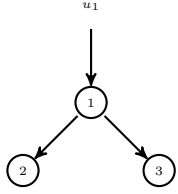
$$A_{1}=
\begin{bmatrix}
  0 & 1 \\
  2 & 0 \\
\end{bmatrix}
,~
A_{2}=
\begin{bmatrix}
  1 & 3 \\
  0 & 1 \\
\end{bmatrix}
,~
A_{3}=
\begin{bmatrix}
  2 & 3 \\
  0 & 2 \\
\end{bmatrix},$$

$$B_{1}=
\begin{bmatrix}
  1 & 0 \\
  0 & 1 \\
\end{bmatrix}
,~
B_{2}=
\begin{bmatrix}
  2 & 0 \\
  0 & 1 \\
\end{bmatrix}
,~
B_{3}=
\begin{bmatrix}
  1 & 0 \\
  0 & 3 \\
\end{bmatrix},$$

$$C_{1}=
\begin{bmatrix}
  0 & 1 \\
  1 & 0 \\
\end{bmatrix}
,~
C_{2}=
\begin{bmatrix}
  2 & 1 \\
  1 & 3 \\
\end{bmatrix}
,~
C_{3}=
\begin{bmatrix}
  3 & 1 \\
  1 & 0 \\
\end{bmatrix},$$

$$H=
\begin{bmatrix}
  0 & 1 \\
  1 & 0 \\
\end{bmatrix},$$
the resulting network is controllable.
\end{examp}

The above two examples reveal some fundamental differences between the controllabilities of heterogeneous networks and homogeneous networks. Therefore, the criteria derived for homogeneous networks might not be applicable to heterogeneous networks. This observation motivates the present study of characterizing necessary and/or sufficient conditions for the controllability of general heterogenous MIMO LTI networks.

\section{Heterogeneous MIMO systems}
In this section, a necessary and sufficient condition is presented for the controllability of the heterogeneous networked MIMO system \eqref{eq:sys_mod_stacked}.
\begin{theorem}\label{theorem-1}
The heterogeneous system \eqref{eq:sys_mod_stacked} is controllable if and only if
\begin{equation}
\left\{ \begin{array}{ll}

\alpha_{i}(sI_{n}-A_{i})-\left( \sum_{j=1,j\ne i}^{N}\omega_{ji}\alpha_{j} \right)HC_{i}=0,\\
\delta_{i}\alpha_{i}B_{i}=0,~~~~i=1,\dots,N.\\

\end{array}\right.\label{eq:theorem-1}
\end{equation}
has a unqiue solution $\alpha_i=0$ for any complex number $s$ and for all $i=1,\dots,N$.

\end{theorem}
\begin{pf}
According to the PBH rank criterion \cite{chen1989}, the system \eqref{eq:sys_mod_stacked} is controllable if and only if the only solution
to
\begin{equation}
\left\{ \begin{array}{ll}
\alpha^{T}\Phi=s\alpha^{T},\\
\alpha^{T}\Psi=0,\\
\end{array}\right. \label{eq:PBH}
\end{equation}
is
$\alpha=0$ for all $s \in \mathbb{C}$. Let $\alpha=[\alpha_{1},\dots,\alpha_{N}]^{T}$ with $\alpha_{i}\in \mathbb{C}^{1 \times n}$. Then, one can rewrite \eqref{eq:PBH} as
\begin{equation}
\left\{\begin{array}{ll}
[\alpha_{1},\dots,\alpha_{N}](sI_{Nn}-\Lambda-\Gamma)=0,\\

[\alpha_{1},\dots,\alpha_{N}]\Psi=0,\\
\end{array}\right.
\end{equation}
which is equivalent to \eqref{eq:theorem-1}. The proof is thus completed.\hfill$\square$
\end{pf}
Theorem \ref{theorem-1} provides a necessary and sufficient condition for the controllability of \eqref{eq:sys_mod_stacked}.
Based on Theorem \ref{theorem-1}, several precise results can be derived for some special networks, which are easily to use.

\begin{cor}\label{cor-1}
Consider a directed chain network of $N$ nodes, with node 1 being the root. Let
$$A_{i}=\begin{bmatrix}
  0 & a_{i1} & \cdots & 0 \\
  0 & 0 & \ddots & \vdots \\
  \vdots & \vdots & \ddots & a_{in-1} \\
  0 & 0 & \cdots & 0 \\
\end{bmatrix},$$
$$HC_{i}=\begin{bmatrix}
  0 & 0 & \cdots & 0 \\
  0 & 0 & \cdots & 0 \\
  \vdots & \vdots & \ddots & \vdots \\
  h_{i1} & 0 & \cdots & 0 \\
\end{bmatrix}, $$
and $B_{1}=[0,0\dots,1]^{T}$. The chain network is controllable if and only if only the root node is under control, i.e., $\delta_1=1$ and $\delta_{2}=\delta_{3}=\dots=\delta_{N}=0$.
\end{cor}

\begin{pf}
    The necessity is obvious, since if the controlled node is not the root then the root hence the chain will not be controllable.

    To show the sufficiency, observe that since the network topology is a directed chain and only the root node is under control, Eq.~ \eqref{eq:theorem-1} can be simplified as
    \begin{equation}
        \left\{ \begin{array}{ll}
            \alpha_{i}(sI_{n}-A_{i})-\omega_{i+1,i}\alpha_{i+1}HC_{i}=0\\
            (i=1,...,N-1),\\
            \alpha_{N}(sI_{n}-A_{N})=0,\\
            \alpha_{1}B_{1}=0.\\
        \end{array}\right.
    \end{equation}
    The equation $\alpha_{1}B_{1}=0$ indicates that the $n$th variable of $\alpha_{1}$ is zero, which, together with $\alpha_{1}(sI_{n}-A_{1})-\omega_{21}\alpha_{2}HC_{1}=0$, yields
    \begin{equation}
        \left\{ \begin{array}{ll}
            \alpha_{1n}=0,\\
            \alpha_{11}s-\omega_{21}\alpha_{2n}h_{11}=0,\\
            \alpha_{12}s-\alpha_{11}a_{11}=0,\\
            \alpha_{13}s-\alpha_{12}a_{12}=0,\\
            ~~~~~~~\vdots\\
            \alpha_{1n}s-\alpha_{1n-1}a_{1n-1}=0,\\
        \end{array}\right.
    \end{equation}
    indicating that $\alpha_{1}=0$ and the $n$th variable $\alpha_{2}$ is zero. Similarly, one obtains $\alpha_{2}=\alpha_{3}=\dots=\alpha_{n-1}=0$, and the $n$th variable $\alpha_{n}$ is zero. From $\alpha_{N}(sI_{n}-A_{N})=0$, it follows that $\alpha_{n}$ is zero. According to Theorem \ref{theorem-1}, the networked system is controllable.\hfill$\square$
%
%
\end{pf}

In what follows, consider the case that there exists one node without incoming edges.


\begin{cor}\label{cor-2}
Suppose that there exists one node $k$ without incoming edges. The heterogeneous networked system \eqref{eq:sys_mod_stacked} is controllable only if $(A_{k},B_{k})$ is controllable, node $k$ is under control, and for any complex number $s$, the only solution of the following equations
    \begin{equation}\label{eq:cor-2}
    \left\{ \begin{array}{ll}

    \alpha_{i}(sI_{n}-A_{i})-\sum_{j=1,j\ne i,j\ne k}^{N}\omega_{ji}\alpha_{j}HC_{i}=0,\\
    \delta_{i}\alpha_{i}B_{i}=0,~~~~i=1,\dots,N, i\neq k,\\

    \end{array}\right.
    \end{equation}
is $\alpha_{i}=0$ for all $i=1,...,N$.
\end{cor}
\begin{pf}
    Since node $k$ does not have any incoming edge, Eq. \eqref{eq:theorem-1} can be rewritten as
    \begin{equation}
        \left\{ \begin{array}{ll}

        \alpha_{i}(sI_{n}-A_{i})-\sum_{j=1,j\ne i,j\ne k}^{N}\omega_{ji}\alpha_{j}HC_{i}=0,\\
        \delta_{i}\alpha_{i}B_{i}=0,~~~~i=1,\dots,N, i\neq k,\\
        \alpha_{k}(sI_{n}-A_{k})=0,\\
        \delta_{k}\alpha_{k}B_{k}=0,\\

        \end{array}\right.
        \end{equation}
    and the $k$th block row of $\Phi$ in (3) becomes $[0,...,A_{k},...,0]$. If $\delta_{k}=0$, then for any $s_{0}\in\sigma(A_{k})$, the row rank of $[s_{0}I_{Nn}-\Phi,\Psi]$ will be reduced at least by one. If $(A_{k},B_{k})$ is uncontrollable, then there exists  $s_{0}\in\sigma(A_{k})$ such that the rank of $(s_{0}I_{n}-A_{k},B_{k})$ is less than $n$, which also reduces the rank of $[s_{0}I_{Nn}-\Phi,\Psi]$. The proof is thus completed.\hfill$\square$
\end{pf}

Corollary \ref{cor-2} presents a necessary condition for determining the controllability of system  \eqref{eq:sys_mod_stacked} when there exists one node without incoming edges. The following theorem gives a necessary condition for generic heterogeneous networks.

\begin{theorem}\label{theorem-2}
Suppose that there exists one node, $i$, without external control inputs. In order to make the heterogeneous networked system \eqref{eq:sys_mod_stacked} controllable,
it is necessary that $[-\omega_{i1}HC_{1},-\omega_{i2}HC_{2},\dots ,sI-A_{i}, \dots ,-\omega_{iN}HC_{N}]$ has full row-rank for any complex number $s$.
\end{theorem}
\begin{pf}
Since node $i$ does not have any external input, the $i$th
block row of $[sI-\Phi,\Psi]$ is $[-\omega_{i1}HC_{1},-\omega_{i2}HC_{2},$ $\dots,sI-A_{i},\dots,-\omega_{iN}HC_{N}]$. If
any row in this block is not independent of the others, it follows that $\mathrm{rank}(sI-\Phi,\Psi)<Nn$, and the heterogeneous networked system is uncontrollable according to the PBH rank criterion \cite{chen1989}.\hfill$\square$
\end{pf}


The effect of the number of external control inputs is explored by the following theorem. Without loss of generality, suppose that the first $m$ nodes are under control.

\begin{theorem}\label{theorem-3}
Suppose $N>\sum_{i=1}^{m}\mathrm{rank}(B_{i})$, $A_{1}, A_{2},\dots, A_{N}$ are similar to each other, and there exist $k_i \neq 0$, $i=1,\dots,N$, such that $k_{1}C_{1}=\dots=k_{N}C_{N}$. For the heterogeneous networked system \eqref{eq:sys_mod_stacked} to be controllable, it is necessary that $(A_{i},C_{i})$ is observable for all i =1,...,N.

\end{theorem}
\begin{pf}
    Assume that there exists a node $i$, such that $(A_{i},C_{i})$ is unobservable. As a result, there must exist a complex number $s_{0}\in\sigma(A_{i})$ and a nonzero vector $\alpha_{i}\in \mathbb{C}^{n}$
    such that
    \begin{equation}
    \left\{\begin{array}{ll}
    C_{i}\alpha_{i}=0,\\
    (s_{0}I_{n}-A_{i})\alpha_{i}=0.\\
    \end{array}\right.\label{eq:obs}
    \end{equation}

    \noindent Let $\Phi_{s_{0}}=[\Phi_{s_{0}}^{1},\dots,\Phi_{s_{0}}^{N}]\triangleq s_{0}I_{Nn}-\Phi$. Then,
        $$\Phi_{s_{0}}^{i}=[\omega_{1i}(HC_{i})^{T},...,(sI_{n}-A_{i})^{T},...,\omega_{Ni}(HC_{i})^{T}]^{T}.$$
    By \eqref{eq:obs}, one has
        $$\Phi_{s_{0}}^{i}\alpha_{i}=0,$$
    which implies that rank$(\Phi_{s_{0}}^{i})\le n-1$. Because $A_{1}, A_{2},\dots, A_{N}$ are similar to each other, they have the same eigenvalue $s_{0}$. Following a similar analysis, one has
    $\mathrm{rank}(\Phi_{s_{0}}^{i})\le n-1, i=1,\dots,N$, i.e., $\mathrm{rank}(\Phi_{s_{0}})\le N(n-1)$, where $k_{1}C_{1}=\dots=k_{N}C_{N}$ is employed. It follows from $\sum_{i=1}^{m}\mathrm{rank}(B_{i})<N$ that $\mathrm{rank}(s_{0}I_{Nn}-\Phi,\Psi)<Nn$,
    which suggests that the heterogeneous networked system \eqref{eq:sys_mod_stacked} is uncontrollable according to the PBH criterion \cite{chen1989}.\hfill$\square$
\end{pf}

The following theorem explores the effect of the network topology on the controllability.

\begin{theorem}\label{theorem-4}
If $A_{1}+s_{0}HC_{1}=...=A_{N}+s_{0}HC_{N}$ for all $s_{0}\in\sigma(W)$, the heterogeneous networked system \eqref{eq:sys_mod_stacked} is controllable only if $(W,\bigtriangleup)$ is
controllable.
\end{theorem}
\begin{pf}
If $(W,\bigtriangleup)$ is uncontrollable, then there exist an $s_{0}\in\sigma(W)$ and a nonzero vector $\xi\in \mathbb{C}^{N\times 1}$ such that
\begin{equation}
\left\{\begin{array}{ll}
\xi^{T}(s_{0}I_{N}-W)=0,\\
\xi^{T}\bigtriangleup=0.\\
\end{array}\right.
\end{equation}
It follows from \eqref{eq:sys_mod_stacked} that
\begin{align}
\begin{split}
(\xi^T \otimes I_{n})\dot x=&(\xi^T \otimes I_{n})((\Lambda+\Gamma)x+\Psi u),
\nonumber
\end{split}
\end{align}
which can be equivalently written as
\begin{align}
\begin{split}
\Bigg(\sum_{i=1}^{N}\xi_{i}x_{i}\Bigg)^{'}=\sum_{i}^{N}A_{i}\xi_{i}x_{i}+\sum_{i}^{N}s_{0}HC_{i}\xi_{i}x_{i},
\nonumber
\end{split}
\end{align}
where $\xi^{T}\bigtriangleup=0$ is used to obtain $(\xi^T \otimes I_n) \Psi=0$, and the superscript $'$ denotes the derivative.
Since $A_{1}+s_{0}HC_{1}=...=A_{N}+s_{0}HC_{N}$, one has
\begin{equation}
\left(\sum_{i=1}^{N}\xi_{i}x_{i}\right)^{'}=(A_{1}+s_{0}HC_{1})\sum_{i}^{N}\xi_{i}x_{i}.
\end{equation}
This implies that the variable $\sum_{i=1}^{N}\xi_{i}x_{i}$ is unaffected by the external control input $u$. For the zero initial state
$x_{i}(t_{0}),i=1,\dots,N$, one has $\sum_{i=1}^{N}\xi_{i}x_{i}(t_{0})=0$. Moreover, it follows from the uniqueness of the solution to the linear equation (11) that $\sum_{i=1}^{N}\xi_{i}x_{i}(t)=0$ for all $t\ge t_{0}$. Consequently, for any state
$\tilde x\triangleq[\tilde x_{1}^{T},...,\tilde x_{N}^{T}]^{T}$ with $\sum_{i=1}^{N}\xi_{i}\tilde x_{i}\ne 0$, there is no external
control input $u$ that can drive the heterogeneous networked system (3) to traverse from the state $0$ to $\tilde x$,  implying that the network is uncontrollable.\hfill$\square$
\end{pf}
Theorem \ref{theorem-4} suggests that only in some special case, the controllability of the network topology is necessary to insure the controllability of heterogeneous networks. It is noteworthy that homogeneous networks are covered by Theorem \ref{theorem-4}.

\section{Networked systems with controllable node systems }

In this section, consider the controllability of heterogeneous networked systems with controllable node systems described by
\begin{equation}\label{eq:one-dim}
\dot x_{i}=A_{i}x_{i}+\sum_{j=1}^{N}\omega_{ij}HCx_{j}+B_{i}u_{i},
\end{equation}
where
$$A_{i}=\begin{bmatrix}
  0 & 1 & 0 & \cdots & 0 \\
  0 & 0 & 1 & \cdots & 0 \\
  \vdots & \vdots & \vdots & \ddots & \vdots \\
  0 & 0 & 0 & \cdots & 1 \\

  -a_{i,0} & -a_{i,1} & -a_{i,2} & \cdots & -a_{i,n-1} \\
\end{bmatrix},
B_{i}=\begin{bmatrix}
  0 \\
  0 \\
  \vdots\\
  0 \\
  1 \\
\end{bmatrix}.$$

Let $u_{i}=a_{i}^{T}x_{i}+\delta_{i}u_{oi}$, where $a_{i}=[a_{i,0}, a_{i,1}, \dots, a_{i,n-1}]^{T}$, $u_{oi}\in \mathbb{R}$ is the external control input, and $\delta_{i}=1$ if node $i$ is under external control and $\delta_{i}=0$ otherwise. Eq. \eqref{eq:one-dim} can be rewritten in terms of the external inputs as
\begin{equation}\label{eq:one-dim2}
\dot x_{i}=Ax_{i}+\sum_{j=1}^{N}\omega_{ij}HCx_{j}+\delta_{i}Bu_{oi},
\end{equation}
where
$$A=A_{i}+B_{i}a_{i}^{T}=\begin{bmatrix}
  0 & 1 & 0 & \cdots & 0 \\
  0 & 0 & 1 & \cdots & 0 \\
  \vdots & \vdots & \vdots & \ddots & \vdots \\
  0 & 0 & 0 & \cdots & 1 \\

  0 & 0 & 0 & \cdots & 0 \\
\end{bmatrix},
B=B_{i}=\begin{bmatrix}
  0 \\
  0 \\
  \vdots\\
  0 \\
  1 \\
\end{bmatrix}.$$

\noindent Denote
$$\Delta=\mathrm{diag}(\delta_{1},\dots,\delta_{N}).$$
The networked system \eqref{eq:one-dim2} can also be rewritten in a compact form as
\begin{equation}\label{eq:one-dim mod stacked}
\dot x=\Phi x+\Psi u_{oi},
\end{equation}
in which $\Phi=I_{n}\otimes A+W\otimes HC$ and $\Psi=\Delta\otimes B$.

\subsection{Networks with one-dimensional communication}

In this section, consider the case that $B\in \mathbb{R}^{n\times 1}$ and $C\in \mathbb{R}^{1\times n}$, i.e., the input and output of the nodes are one-dimensional. Let the set of nodes under external control be
$$
\nu=\{i=1,...,m ~|~ \delta_{i}\neq0\}, 1\leq m \leq N.
$$
\noindent For $s\in\sigma(A_{i}+B_{i}a_{i}^{T})$, define a matrix set

\[ \alpha(s)= \Set{[\alpha_{1},...,\alpha_{N}]}
    {
        \begin{lgathered}
            \alpha_{i}\in \alpha^{1}(s) ~\mathrm{for} ~i\notin\nu\\
            \alpha_{i}\in \alpha^{2}(s) ~\mathrm{for} ~i\in \nu
        \end{lgathered}
    }
\]
where
\begin{align}
\begin{split}
\alpha^{1}(s) &=\Big\{\xi\in \mathbb{C}^{n\times 1} ~|~ \xi^{T} (sI_{n}-A_{i}-B_{i}a_{i}^T)=0\Big\}\\
\alpha^{2}(s) &=\Big\{\xi\in \mathbb{C}^{n\times 1} ~|~ \xi^{T} B_{i}=0, \xi\in\alpha^{1}(s)\Big\}.
\nonumber
\end{split}
\end{align}

\begin{theorem}\label{theorem-5}
Suppose that $|\nu|<N$. With control input $u_{i}=a_{i}^{T}x_{i}+\delta_{i}u_{oi}$, the networked system \eqref{eq:one-dim} is controllable if and only if the following conditions hold:
\begin{enumerate}
    \item[(i)] $(A_{i}+B_{i}a_{i}^T,H)$ is controllable;

    \item[(ii)] $(A_{i}+B_{i}a_{i}^T,C)$ is observable;

    \item[(iii)] for any $s\in\sigma(A_{i}+B_{i}a_{i}^T)$ and $\rho\in\alpha(s)$, $W^{T}\rho^{T}\neq0$ if $\rho\neq0$;

    \item[(iv)] for any $s\notin\sigma(A_{i}+B_{i}a_{i}^T)$, $\mathrm{rank}(I-W\gamma,\Delta\eta)=N$, where $\gamma=C(sI_{n}-A_{i}-B_{i}a_{i}^T)^{-1}H$ and $\eta=C(sI_{n}-A_{i}-B_{i}a_{i}^T)^{-1}B_{i}$,
\end{enumerate}
where $i=1,...,N$.
\end{theorem}

The following two results from \cite{wang2016} will be employed in the proof of Theorem \ref{theorem-5}.

\begin{lemma}\cite{wang2016}\label{lemma:ah1}
 If there exists one node without external control inputs, then for networked system \eqref{eq:one-dim mod stacked} to be controllable, it is necessary that $(A,HC)$ is controllable.

\end{lemma}

\begin{lemma}\cite{wang2016}\label{lemma:ac}
If the number of nodes with external control inputs is m, and  $N>m\cdot \mathrm{rank}(B)$, then for the networked system \eqref{eq:one-dim mod stacked} to be controllable, it is necessary that $(A,C)$ is observable.
\end{lemma}

Next, a new result is established.

\begin{lemma}\label{lemma:ah2}
Suppose that $C\in \mathbb{R}^{1\times n}$ and $H\in \mathbb{R}^{n\times 1}$ are non-zero. Then, $(A,HC)$ is controllable if and only if $(A,H)$ is controllable.
\end{lemma}
\begin{pf}
Since $\mathrm{rank}(H)+\mathrm{rank}(C)-1\leqslant \mathrm{rank}(HC)\leqslant \min\{\mathrm{rank}(H),\mathrm{rank}(C)\}$, one has $\mathrm{rank}(HC)=1$. Therefore, $\mathrm{rank}(sI-A,HC)=\mathrm{rank}(sI-A,H)$, which leads to the conclusion.\hfill$\square$
\end{pf}

Now, it is ready to prove Theorem \ref{theorem-5}.

\begin{pf}
(Necessity): The assumption $|\nu|<N$ indicates that there exists at least one node without external control inputs. It follows from Lemmas \ref{lemma:ah1} and \ref{lemma:ah2} that the controllability of $(A_{i}+B_{i}a_{i}^{T},H)$ is necessary for the controllability of system \eqref{eq:one-dim mod stacked}, which shows the necessity of condition $(i)$.  The condition (ii) follows readily from Lemma \ref{lemma:ac}.

Now, suppose that the condition (iii) is not satisfied. Then, there exist an $s_{0}\in \sigma(A _{i}+B_{i}a_{i}^{T})$ and  a non-zero matrix $\rho\in \alpha(s_{0})$ such that
$$W^{T}\rho^{T}=0,$$
 in which $\rho=[\rho_{1}^{T},\dots,\rho_{N}^{T}]$ and $\rho_{i}\in \mathbb{R}^{1\times n}$. The equality $W^{T}\rho^{T}=0$ is equivalent to $\sum_{j=1,j\neq i}^{N}\omega_{ji}\rho_{j}=0$, $i=1,\dots,N$. Let $\bar{\rho}=[\rho_{1},\dots,\rho_{N}]$. It is straightforward to verify that $\bar{\rho}\Psi=0$ and
\begin{equation}\begin{split}
& \quad \bar{\rho}(s_{0}I_{Nn}-\Phi)\\
                &=\bar{\rho}(s_{0}I_{Nn}-I_{N}\otimes A-W\otimes HC)\\
                 &=\bar{\rho}(s_{0}I_{Nn}-I_{N}\otimes (A_{i}+B_{i}a_{i}^T)-W\otimes HC)\\
                 &=-\bar{\rho}(W\otimes HC)\\
                 &=-\left[\left(\sum_{j=1,j\neq 1}^{N}\omega_{j,1}\rho_{j}\right)HC,\dots, \left(\sum_{j=1,j\neq N}^{N}\omega_{j,N}\rho_{j} \right) HC\right] \\
                 &=0, \nonumber
\end{split}
\end{equation}
which implies that the network is uncontrollable.

Finally, suppose that the condition (iv) is not satisfied. Then, there exists an $s_{0}\notin\sigma(A_{i}+B_{i}a_{i}^T)$ satisfying
$$\mathrm{rank}(I_{N}-W\gamma_{0},\Delta\eta_{0})<N,$$
with $\gamma_{0}=C(s_{0}I_{n}-A_{i}-B_{i}a_{i}^T)^{-1}H$ and $\eta_{0}=C(s_{0}I_{n}-A_{i}-B_{i}a_{i}^T)^{-1}B_{i}$.
Thus, there exists a non-zero vector $\xi=[\xi_{1},...,\xi_{N}]\in \mathbb{C}^{1\times N}$ such that
\begin{align}
\xi(I_{N}-W\gamma_{0})=0~~~\text{and}~~~\xi\Delta\eta_{0}=0. \label{eq:b1}
\end{align}
Let $\alpha=[\alpha_{1},...,\alpha_{N}]$ with $\alpha_{i}=\xi_{i}C(s_{0}I-A_{i}-B_{i}a_{i}^T)^{-1}$. Since $\xi\neq 0$, one has $\alpha\neq 0$. It follows from \eqref{eq:b1} that
\begin{equation}\begin{split}
\alpha\Psi&=(\xi\otimes C(s_{0}I_{n}-A_{i}-B_{i}a_{i}^T)^{-1})(\Delta\otimes B)\\
          &=(\xi\Delta)\otimes(C(s_{0}I_{n}-A_{i}-B_{i}a_{i}^T)^{-1}B_{i})\\
          &=\xi\Delta\eta_{0}\\
          &=0 \nonumber
\end{split}
\end{equation}
and
\begin{equation}\begin{split}
		& \quad \alpha(s_{0}I_{Nn}-\Phi) \\
        &=(\xi\otimes C(s_{0}I_{n}-A_{i}-B_{i}a_{i}^T)^{-1})\\
        & \quad \times (I_{N}\otimes(s_{0}I_{n}-A)-W\otimes HC)\\
        &=\xi\otimes C-\xi W\otimes(C(s_{0}I_{n}-A)^{-1}H)C\\
        &=\xi(I-W\gamma_{0})\otimes C\\
        &=0, \nonumber
\end{split}
\end{equation}
which implies that the network is uncontrollable.

\noindent (Sufficiency): For $s\in \mathbb{C}$, suppose that there exists a vector $\alpha=[\alpha_{1},...,\alpha_{N}]$ with $\alpha_{i}\in \mathbb{C}^{1\times n}$ such that $\alpha(sI-\Phi)=0$ and $\alpha\Psi=0$, which are equivalent to
\begin{equation}\label{eq:suff-sys1}
\alpha_{i}(sI_{n}-A)-\sum_{j=1,j\ne i}^{N}\omega_{ji}\alpha_{j}HC=0,~~~i=1,...,N,
\end{equation}
and
\begin{equation}\label{eq:suff-sys2}
\alpha_{i}B=0,~~~~i\in\nu.
\end{equation}
If $s\in\sigma(A_{i}+B_{i}a_{i}^T)$, it follows that $\mathrm{rank}(sI_{n}-A_{i}-B_{i}a_{i}^T)=\mathrm{rank}(sI_{n}-A)<n$. Eq.~\eqref{eq:suff-sys1} yields that, for all $i=1,...,N$,
\begin{equation}\label{eq:suff-con1}
\sum_{j=1,j\neq i}^{N}\omega_{ji}\alpha_{j}H=0.
\end{equation}
Otherwise, $C$ is a linear combination of the row vectors of $sI_{n}-A$, yielding
$$\mathrm{rank}\begin{pmatrix}
\begin{bmatrix}
  C \\
  sI_{n}-A\\
\end{bmatrix}
\end{pmatrix}=\mathrm{rank}(sI_{n}-A)<n.
$$
Since $A=A_{i}+B_{i}a_{i}^T$, it contradicts the observability of $(A_{i}+B_{i}a_{i}^T,C)$. Substituting \eqref{eq:suff-con1} into \eqref{eq:suff-sys1} leads to
\begin{equation}\label{eq:suff-con2}
\alpha_{i}(sI_{n}-A)=0.
\end{equation}
Therefore, for all $i=1,...,N$, one has
\begin{equation}
\sum_{j=1,j\neq i}^{N}\omega_{ji}\alpha_{j}(sI_{n}-A)=0.
\end{equation}
Combining \eqref{eq:suff-con1} and the controllability of $(A_{i}+B_{i}a_{i}^T,H)$ gives
\begin{equation}\label{eq:suff-con3}
\sum_{j=1,j\neq i}^{N}\omega_{ji}\alpha_{j}=0.
\end{equation}

Next, let $\rho=[\alpha_{1}^{T},...,\alpha_{N}^{T}]$. In light of \eqref{eq:suff-sys2}, \eqref{eq:suff-con2} and \eqref{eq:suff-con3}, it can be verified that $W^{T}\rho^{T}=0$ with $\alpha_{i}(sI_{n}-A)=0, i=1,...,N$, and $\alpha_{i}B=0, i\in\nu$. Therefore, by condition (iii), one has $\alpha_{i}=0$ for all $i=1,...,N$.

If $s\notin\sigma(A_{i}+B_{i}a_{i}^T)$, then $(sI_{n}-A)=(sI_{n}-A_{i}-B_{i}a_{i}^T)$ is invertible. It follows  from \eqref{eq:suff-sys1} that
\begin{equation}
\alpha_{i}=\sum_{j=1,j\neq i}^{N}\omega_{ji}\alpha_{j}HC(sI_{n}-A)^{-1},~~~i=1,...,N.
\end{equation}
Let $\xi_{i}=\sum_{j=1,j\neq i}^{N}\omega_{ji}\alpha_{j}H$. Then, for $i=1,...,N$, one has
\begin{equation}\label{eq:suff-alphaxi}
\alpha_{i}=\xi_{i}C(sI_{n}-A)^{-1}
\end{equation}
and
\begin{equation}\label{eq:suff-con4}
\begin{split}
\xi_{i} &=\sum_{j=1,j\neq i}^{N}\omega_{ji}\alpha_{j}H\\
        &=\sum_{j=1,j\neq i}^{N}\omega_{ji}\xi_{j}C(sI_{n}-A)^{-1}H\\
        &=\sum_{j=1,j\neq i}^{N}\omega_{ji}\xi_{j}\gamma,
\end{split}
\end{equation}
where $\gamma$ is defined in (iv) of the theorem. Let $\xi=[\xi_{1},...,\xi_{N}]$, and rewrite \eqref{eq:suff-con4} as
\begin{equation}\label{eq:suff-xileft}
\xi(I_{n}-W\gamma)=0.
\end{equation}
It follows from \eqref{eq:suff-sys2} and \eqref{eq:suff-alphaxi} that $\xi_{i}C(sI_{n}-A)^{-1}B=0$ for $i\in\nu$, which is equivalent to
\begin{equation}
\xi\Delta\eta=0,
\end{equation}
where $\eta$ is defined in (iv) of the theorem. 
Combining \eqref{eq:suff-xileft} and condition (iv) leads to $\xi=0$, which implies that $\alpha_{i}=0$ for all $i=1,...,N$, by \eqref{eq:suff-alphaxi}.

From the above analysis, for any $s\in \mathbb{C}$, the row vectors of the matrix $[sI_{Nn}-\Phi,\Psi]$ are linearly independent, implying that rank$(sI_{Nn}-\Phi,\Psi)=N n$. Thus, the networked system \eqref{eq:one-dim} is controllable.\hfill$\square$
\end{pf}

The following provides a sufficient condition for the controllability of \eqref{eq:one-dim}.
In order to state the main results, define
\[ \kappa= \Set{[\underbrace{0,...,0}_{n-1},\kappa_{n}]^{T}}
{
    \kappa_{n}\in \mathbb{R}, \kappa_{n}\neq0
}
\]
\[ \tau= \Set{[\tau_{1},\underbrace{0,...,0}_{n-1}]}
{
    \tau_1\in \mathbb{R}, \tau_1\neq0
}.
\]
\begin{cor}\label{cor-3}
Suppose that $|\nu|<N$ and $u_{i}=a_{i}^{T}x_{i}+\delta_{i}u_{oi}$. The networked system \eqref{eq:one-dim} is controllable if the following three conditions hold:

\begin{enumerate}
    \item[(i)] $H\in\kappa$ and $C\in\tau$;

    \item[(ii)] for $\rho\in\alpha(0)$, $W^{T}\rho^{T}\neq0$ if $\rho\neq0$;

    \item[(iii)] for any $s\neq 0$, $\mathrm{rank}(I-W\gamma,\Delta\eta)=N$, where $\gamma=C(sI_{n}-A_{i}-B_{i}a_{i}^T)^{-1}H$ and $\eta=C(sI_{n}-A_{i}-B_{i}a_{i}^T)^{-1}B_{i}$,
\end{enumerate}
where $i=1,...,N$.
\end{cor}

\begin{pf}
Assume that $\xi=[\xi_{1},\dots,\xi_{n}]\in \mathbb{R}^{1\times n}$. For any $s\in \mathbb{C}$, let $\xi(sI_{n}-A_{i}-B_{i}a_{i}^T)=0$ and $\xi H=0$. It follows that
\begin{equation}\label{eq:cor-3}
    \left\{ \begin{array}{ll}

    s\xi_{1} =0,\\
        s\xi_{2}-\xi_{1} =0,\\
        s\xi_{3}-\xi_{2} =0,\\
        ~~~~~\vdots\\
        s\xi_{n}-\xi_{n-1} =0,\\
        \xi_{n}\kappa_{n} =0.
        \nonumber

    \end{array}\right.
\end{equation}
If $s\neq 0$, one has $\xi_{1}=,\dots,=\xi_{n}=0$ and $\xi_{n}\kappa_{n}=0$, which yields $\xi=0$. If $s=0$, one has $\xi_{1}=,\dots,=\xi_{n-1}=0$ and $\xi_{n}\kappa_{n}=0$, which also leads to $\xi=0$ since $\kappa_{n}\neq 0$. Therefore, $(A_{i}+B_{i}a_{i}^T,H)$ is controllable. Similarly, it can be verified that $(A_{i}+B_{i}a_{i}^T,C)$  is observable. In addition, the determinant of the matrix $(sI_{n}-A_{i}-B_{i}a_{i}^T)$ is
$$
\mathrm{det}(sI_{n}-A_{i}-B_{i}a_{i}^T)=\mathrm{det}\left(\begin{bmatrix}
  s & -1 & 0 & \cdots & 0 \\
  0 & s & -1 & \cdots & 0 \\
  \vdots & \vdots & \vdots & \ddots & \vdots \\
  0 & 0 & 0 & \cdots & -1 \\
  0 & 0 & 0 & \cdots & s \\
\end{bmatrix}
\right)\\
.$$

It can be shown that $\mathrm{det}(sI_{n}-A_{i}-B_{i}a_{i}^T)=s^{n}$ by using the expansion rule for calculating matrix determinant. It thus follows that all the eigenvalues of $(A_{i}+B_{i}a_{i}^T)$ are zeros, which implies that conditions (ii) and (iii) in Corollary \ref{cor-3} are equivalent to conditions (iii) and (iv) in Theorem \ref{theorem-5}. Thus, the four conditions of Theorem \ref{theorem-5} are satisfied by the conditions of Corollary \ref{cor-3}.\hfill$\square$
\end{pf}

\subsection{Networks with diagonalizable network topologies}

In this section, consider the controllability of the networked system \eqref{eq:one-dim}, in which $W$ is assumed to be diagonalizable. First, recall a lemma from \cite{hao2017}.

\begin{lemma}\cite{hao2017}\label{lemma:diagW}
Assume that $W$ is diagonalizable with eigenvalues $\lambda_{1},\dots,\lambda_{N}$. Denote $\Theta=\{\lambda_{1},\dots,\lambda_{N}\}$. The networked system \eqref{eq:one-dim} is controllable if and only if

\begin{enumerate}
    \item[(i)] $(W,\Delta)$ is controllable;

    \item[(ii)] $(A+\lambda_{i}H,B)$ is controllable, for $i=1,2,\dots,N$;

    \item[(iii)] if matrices $A+\lambda_{i_{1}}H,\dots,A+\lambda_{i_{p}}H(\lambda_{i_{k}}\in\Theta, for~ k = 1,\dots,p,p > 1)$ have a common eigenvalue $\theta$, then $\underbrace{(t_{i_{1}}\Delta)\otimes(\xi_{i_{1}}^{1}B),\dots,(t_{i_{1}}D)\otimes(\xi_{i_{1}}^{\gamma_{i_{1}}})}_{\text{Geometric multiplicity is}~ \gamma_{i_{1}}},\dots,$\\
        $\underbrace{(t_{i_{p}}\Delta)\otimes(\xi_{i_{p}}^{1}B),\dots,(t_{i_{p}}D)\otimes(\xi_{i_{p}}^{\gamma_{i_{p}}})}_{\text{Geometric multiplicity is}~ \gamma_{i_{p}}}$ are linearly independent, where $t_{i_{k}}$ is the left eigenvector of $W$ corresponding to the eigenvalue $\lambda_{i_{k}}$; $\gamma_{i_{k}}\geqslant 1$ is the geometric multiplicity of $\theta$ for$A+\lambda_{i_{k}}H$; $\xi_{i_{k}}^{l}(l = 1,\dots,\gamma_{i_{k}})$ are the left eigenvectors of $A+\lambda_{i_{k}}H$ corresponding to $\theta,k = 1,2,\dots,p$.
\end{enumerate}

\end{lemma}

The following theorem determines the controllability of the networked system \eqref{eq:one-dim} with a diagonalizable $W$.

\begin{theorem}
Assume that $W$ is diagonalizable with eigenvalues $\Theta=\{\lambda_{1},\dots,\lambda_{N}\}$. Under control inputs $u_{i}=a_{i}^{T}x_{i}+\delta_{i}u_{oi}$, the networked system \eqref{eq:one-dim} is controllable if and only if

\begin{enumerate}
    \item[(i)] $(W,\Delta)$ is controllable;

    \item[(ii)] $(A_{i}+B_{i}a_{i}^T+\lambda_{i}H,B_{i})$ is controllable, for $i=1,2,\dots,N$;

    \item[(iii)] if matrices $A_{i_{1}}+B_{i_{1}}a_{i_{1}}+\lambda_{i_{1}}H,\dots,A_{i_{p}}+B_{i_{p}}a_{i_{p}}+\lambda_{i_{p}}H(\lambda_{i_{k}}\in\Theta, \mathrm{for}~k = 1,\dots,p,p > 1)$ have a common eigenvalue $\theta$, then \\ $\underbrace{(t_{i_{1}}\Delta)\otimes(\xi_{i_{1}}^{1}B),\dots,(t_{i_{1}}D)
        \otimes(\xi_{i_{1}}^{\gamma_{i_{1}}})}_{\text{geometric multiplicity is}~ \gamma_{i_{1}}},$
        $\dots,$ \\ $\underbrace{(t_{i_{p}}\Delta)\otimes(\xi_{i_{p}}^{1}B),\dots,(t_{i_{p}}D)\otimes(\xi_{i_{p}}^{\gamma_{i_{p}}})}_{\text{geometric multiplicity is}~ \gamma_{i_{p}}}$ are linearly independent, where $t_{i_{k}}$ is the left eigenvector of $W$ corresponding to the eigenvalue $\lambda_{i_{k}}$; $\gamma_{i_{k}}\geqslant 1$ is the geometric multiplicity of $\theta$ for$A+\lambda_{i_{k}}H$; $\xi_{i_{k}}^{l}(l = 1,\dots,\gamma_{i_{k}})$ are the left eigenvectors of $A+\lambda_{i_{k}}H$ corresponding to $\theta,k = 1,2,\dots,p$.
\end{enumerate}

\end{theorem}

\begin{pf}
With control inputs $u_{i}=a_{i}^{T}x_{i}+\delta_{i}u_{oi}$, the controllability of the system \eqref{eq:one-dim} is equivalent to the controllability of the system \eqref{eq:one-dim mod stacked}. The result follows immediately from Lemma~\ref{lemma:diagW}.\hfill$\square$
\end{pf}

\section{Conclusion}
This study has revealed how the network topology, the external control input, the node system $(A_{i},B_{i},C_{i})$ as well as the inner interaction matrix $H$ affect the controllability of directed and weighted heterogeneous networked MIMO LTI systems. A necessary and sufficient condition has been derived for determining whether such a  network is controllable. The results demonstrate that the heterogeneity of node dynamics is a fundamental factor for the network controllability. It is found that, in some specific cases,  the controllability of $(A_{i},B_{i})$ is necessary for the controllability of the heterogeneous network. Another interesting finding is that the heterogeneous networked system can be controllable even if the network topology is uncontrollable. This work might lead to a better understanding or even manipulation of the controllability of complex networked systems, especially in the
heterogeneous setting.

%

\end{document}